\documentclass[a4paper,twoside]{amsart}

\usepackage{a4wide,amscd}

\newtheorem{theorem}{Theorem}
\newtheorem{corollary}{Corollary}
\newtheorem{lemma}{Lemma}
\newtheorem{proposition}{Proposition}
\theoremstyle{definition}
\newtheorem{definition}{Definition}

\theoremstyle{remark}

\def\htop{h_{\mathrm{top}}}
\def\hhom{h_{\mathrm{hom}}}
\def\hper{h_{\mathrm{per}}}

\DeclareMathOperator{\tr}{tr}
\DeclareMathOperator{\Fix}{Fix}
\DeclareMathOperator{\Image}{Im}
\DeclareMathOperator{\Ker}{Ker}
\DeclareMathOperator{\Ver}{Ver}

\title{Topological entropy, homological growth and zeta functions on graphs}

\author[J.~F.~Alves, R.~Hric and J.~Sousa Ramos]{Jo\~ao Ferreira Alves$^1$,
Roman Hric$^{1,2}$ and Jos\'e Sousa Ramos$^1$}

\address{$^1$Center for Mathematical Analysis, Geometry, and Dynamical Systems, 
Instituto Superior T\'ecnico,
Av. Rovisco Pais,
P--1049-001 Lisboa, Portugal}
\email{jalves@math.ist.utl.pt, rhric@math.ist.utl.pt, sramos@math.ist.utl.pt}
\address{$^2$Institute of Mathematics and Computer Science,
Faculty of Natural Sciences,
Matej Bel University,
Tajovsk\'eho 40,
SK--974 01 Bansk\'a Bystrica, Slovakia}
\email{hric@fpv.umb.sk}

\thanks{The first and third authors have been partially supported by FCT/POCTI/FEDER (Portugal).
The second author has been supported by the FCT/POCTI/FEDER, the project SFRH/BPD/10157/2002
(Portugal) and by the VEGA grant 1/0265/03 (Slovakia).}

\begin{document}

\begin{abstract}
In connection with the Entropy Conjecture it is known that the topological entropy of
a continuous graph map is bounded from below by the spectral radius of the induced
map on the first homology group. We show that in the case of a piecewise monotone
graph map, its topological entropy is equal precisely to the maximum of the mentioned
spectral radius and the exponential growth rate of the number of periodic points of negative type.
This nontrivially extends a result of Milnor and Thurston
on piecewise monotone interval maps. For this purpose we generalize the concept
of Milnor-Thurston zeta function incorporating in the Lefschetz zeta function.
The methods developed in the paper can be used also in a more general setting.
\end{abstract}

\maketitle

\section{Introduction and the main results}

One of the most exciting problems in the theory of dynamical systems for the last three decades
has been so called the Entropy Conjecture which was first stated by M. Shub in \cite{S} and
which claims that for a smooth map on a compact differentiable manifold, its topological
entropy is bounded from below by the logarithm of the spectral radius of the induced map in
the corresponding full homology group (called also homological entropy).
Many special cases of this conjecture have been already proved even in the nonsmooth case
but it still remains open in general. In particular, Manning in \cite{M} proved that
the conjecture holds for all continuous maps
on differentiable manifolds if we reduce our attention from the full homology group to
the first homology group. Using the arguments from the paper it is possible to extend
this result for more general spaces including graphs (see \cite{FM}). In this paper
we show that for piecewise monotone graph maps, the topological entropy is not only bounded
from below by its homological entropy but it is exactly equal to the maximum of its homological
entropy and entropy given by the growth of its periodic points of negative type.
(In fact, using this result it is possible to provide an alternative proof of
the Manning's result for continuous graph maps.)

First we recall some notions and definitions needed in the sequel.
By the {\sf\em homological entropy} of a map $f: X \to X$ we mean a topological invariant
$\hhom(f)$ coming from considering the induced linear maps
$f_{\ast i}$ on the homology groups $H_i(X,\mathbb{R})$ and defined by 
\begin{equation*}
\hhom(f)=\log r(f)
\end{equation*}
where $r(f)=\max \{ r(f_{\ast i}): i=0,\dots,\dim X\}$
and $r(f_{\ast i})$ denotes the spectral radius of $f_{\ast i}$.
The purpose of this paper is to establish a precise relationship between
the {\sf\em topological entropy} $\htop(f)$ and the homological entropy $\hhom(f)$ for a
piecewise monotone graph map $f$ (for the definition of topological entropy see any standard
textbook on dynamical systems; a nice introduction into the topic can be found in
\cite{ALM}).
The compact interval and the circle are the simplest examples of graphs.
In general, a {\sf\em graph} is a compact Hausdorff space which
can be written as a union of finitely many homeomorphic copies of the closed
interval $[0,1]$ any two of which intersect at most at their endpoints.
A point of a graph is called its {\sf\em vertex}
if it does not have any open neighborhood homeomorphic to the open interval $]0,1[$.
The set of all vertices of $G$ is denoted by $\Ver(G)$.
Notice that if $G$ is a  graph and $f:G \to G$ is a continuous map, since $H_i(G,\mathbb R)=0$
for $i \geq 2$, and $r(f_{\ast 0})=1$, we obtain
\begin{equation*}
r(f) = \max \{1,r(f_{\ast 1})\}.
\end{equation*}

\begin{definition}
Let $G$ be a graph. A continuous map $f:G \to G$ is called a
{\sf\em piecewise monotone graph} (shortly {\sf\em PMG}) map
if there is a finite set $C \subseteq G$
such that $f$ is injective on each connected component of $G \setminus C$.
\end{definition}

As mentioned before, our goal is to study the relationship between 
$\htop(f)$ and $\hhom(f)$. To this end we define another
topological invariant $\hper^-(f)$. Let $f:G \to G$ be a PMG map. By $\Fix(f)$ we denote
the set of all fixed points of $f$. A point $x \in \Fix(f) \setminus \Ver(G)$
is called {\sf\em of negative type} if $f$ reverses orientation throughout a small neighborhood of
$x$. (Since $x \notin \Ver(G)$, it has a neighborhood homeomorphic to an open real interval
on which we can consider $f$ to be a selfmap of the real line.)
We denote by $\Fix^-(f)$ the set of all fixed points of negative type of $f$. 
Evidently, the set $\Fix(f)$ may be infinite but, since $f$ is a PMG map, the set $\Fix^-(f)$
is always finite. Notice that every iterate of a PMG map
is again a PMG map. Hence the sets $\Fix^-(f^n)$ are always finite 
and therefore we can introduce another topological invariant
\begin{equation*}
\hper^- = \limsup_{n \to \infty} \frac1n \log^+ \#\Fix^-(f^n),
\end{equation*}
the exponential growth rate of the number of periodic points of negative type (we put
$\log^+ x = \log \max \{1,x\}$).

Milnor and Thurston showed in \cite{MT} (see also Theorem 4.11 of \cite{MTr}) that 
\begin{equation}\label{eq002}
\htop(f)=\hper^-(f)
\end{equation}
for any piecewise monotone interval map $f$. Among graph maps this does not hold
anymore --- as an example consider the circle $S^1 = \{ x \in \mathbb C: |x| =1\}$
and the map $f: S^1 \to S^1$ defined by $f(x) = x^2$. For this map one obtains
$\htop(f) = \log 2$ and $\hper^-(f) = 0$. Nevertheless, we prove the following nice relation
extending the last equality.

\begin{theorem}\label{t1}
Let $f$ be a PMG map. Then 
\begin{equation*}
\htop(f) = \max \{ \hper^-(f), \hhom(f) \}.
\end{equation*}
\end{theorem}

The spectral radius $r(f)$ is an algebraic number for any PMG map $f$. Using this we get
as a consequence of the last theorem the next result showing that
the following entropies are equal for almost all values of topological entropy.

\begin{corollary}
Let $f$ be a PMG map and suppose that $\exp(\htop(f))$ is a transcencental number.
Then
\begin{equation*}
\htop(f) = \hper^-(f).
\end{equation*}
\end{corollary}

Notice that the both results hold for PMG maps in general, even for those with
$\Fix(f^n)$ infinite.
In the case that $\Fix(f^n)$ is finite for every $n \geq 1$ then we can consider
 a topological invariant
\begin{equation*}
\hper(f) = \limsup_{n \to \infty} \frac1n \log^+ \#\Fix(f^n).
\end{equation*}
For many important cases of PMG maps (expanding maps and more generally maps
with ``few'' stable periodic orbits), topological entropy represents the
exponential growth rate for the number of periodic orbits, that is
$\htop(f)=\hper(f)$. Moreover, if $\exp(\htop(f))$ is transcendental, we get from the last corollary
the following relation
\begin{equation}
\htop(f)=\hper^-(f)=\hper(f).
\end{equation}
Just stated identity shows that topological entropy in some sense
describes the periodic structure of the system in both quantitative and
qualitative ways --- for an expanding piecewise monotone interval map, we have an obvious
relationship between the number of fixed points of negative and positive types
(the latter one defined analogously)
because between any two consecutive fixed
points of $f^n$ of negative type there is exactly one of its fixed point of positive type
and consequently $\hper^-(f)=\hper(f)$. Indeed, we have no such relation
between the number of  the fixed points of a PMG map of negative and positive types
even if the map is expanding.

One of extremely useful tools for studying the relation between topological entropy
and the growth of the number of periodic points was
introduced by Artin and Mazur in \cite{AM}. Let $X$ be an arbitrary set and $f: X \to X$.
The {\sf\em orbit} of a point $x \in X$ under
the action of $f$ is defined as the set $o_x = \{f^n(x): n \geq 0\}$.
An orbit $o_x$ is said to be {\sf\em periodic} if there is a positive
integer $n$ such that $f^n(x)=x$; the smallest such number we denote by
$\mathrm{p}(o_x)$ and call its {\sf\em period}.
The set of all periodic orbits of $f$ is denoted by $O$. 
Suppose that each positive iterate $f^n$ has
only finitely many fixed points. Then we define the {\sf\em Artin-Mazur zeta
function} of $f$, $\zeta$, to be the formal power series 
\begin{equation*}
\zeta(z) = \exp \sum_{n \geq 1} \frac{\#\Fix(f^n)}n z^n.
\end{equation*}
Recall that the Artin-Mazur zeta function of $f$ is a convenient way of
enumerating the periodic orbits of $f$. Indeed, if each positive iterate of $f$
has only finitely many fixed points then the subset $\{o \in O: \mathrm{p}(o)=k\}$
is for any $k$ always finite and the identity 
\begin{equation*}
\zeta(z) = \prod_{o \in O}\left(1-z^{\mathrm{p}(o)}\right)^{-1}
\end{equation*}
holds in $\mathbb{Z}[[z]]$, the ring of all formal power series in $z$ over
$\mathbb Z$.

Later on, several variants of this notion were introduced by different
authors (cf. \cite{MT}, \cite{BR}; for an extensive survey of the
topic see \cite{Ba}, \cite{P}; cf. also \cite{R}).
In particular, Milnor and Thurston in \cite{MT} modified
the Artin-Mazur zeta function to obtain more information for a piecewise
monotone interval map. Let us very briefly remind how they arrived at the
identity (\ref{eq002}). If $f: I \to I$ is a
piecewise monotone interval map, we call the formal power series 
\begin{equation}\label{eq001}
\zeta^{MT}(z) = \exp \sum_{n \geq 1} \frac{2\#\Fix^-(f^n)-1}n z^n
\end{equation}
the {\sf\em Milnor-Thurston zeta function} of an interval map $f$. Denote its radius of convergence by $\rho$.
Starting from the main relation between $\zeta^{MT}(z)$ and the kneading determinant of $f$,
Milnor and Thurston proved that
\begin{equation*}
\htop(f) = -\log \rho = \hper^-(f).
\end{equation*}

Here we follow the same strategy to prove Theorem \ref{t1}.
As the first step we generalize the concept of Milnor-Thurston zeta function.
Let us begin by defining the Lefschetz and negative zeta functions of a PMG map.
Let $f:G \to G$ be a PMG map. Recall that the formal power series
\begin{equation*}
\zeta^L(z) = \exp \sum_{n\geq 1}\frac{\tr(f_{\ast 0})^n -
\tr(f_{\ast 1})^n}n z^n
\end{equation*}
is called the {\sf\em Lefschetz zeta function} of $f$. We define the
{\sf\em negative zeta function} of $f$ as
\begin{equation*}
\zeta^-(z) = \exp \sum_{n \geq 1} \frac{2\#\Fix^-(f^n)}n z^n.
\end{equation*}
Observe that if $f:I \to I$ is a piecewise monotone interval map then we have
$\tr(f_{\ast 0})=1$ and $\tr(f_{\ast 1}) = 0$ for all $n \geq 1$ and therefore
\begin{equation*}
\zeta^-(z) \zeta^L(z)^{-1} = \exp \sum_{n \geq 1} \frac{2\#\Fix^-(f^n)-1}n z^n
\end{equation*}
holds in $\mathbb Z[[z]]$. So, according to (\ref{eq001}) it is natural to define
the {\sf\em Milnor-Thurston zeta function} of a PMG map $f$ as the formal power series
\begin{equation}\label{eq1}
\zeta^{MT}(z) = \zeta^-(z) \zeta^L(z)^{-1}
\end{equation}
and, as before, there is a close relation between $\htop(f)$ and the radius of convergence
of $\zeta^{MT}(z)$. Theorem \ref{t1} is then an immediate consequence of the following

\begin{theorem}\label{t2}
Let $f$ be a PMG map and denote by $\rho$ the radius of convergence of $\zeta^{MT}(z)$.
Then $0 < \rho \leq 1$ and
\begin{equation*}
\htop(f) = -\log \rho.
\end{equation*}
\end{theorem}

\section{Proof of Theorem \ref{t2}}

The rest of the paper is devoted to the proof of Theorem \ref{t2}. In
order to simplify notation it is convenient to regard a PMG map $f$
as a real map $F$ with discontinuities defined on a subset
of the real line. The proof of Theorem \ref{t2} is
given in two main steps. The first one is the construction
of kneading determinant, $D(z)$, associated to the map $F$. In the second
one we set up the relationship between the kneading determinant and the
zeta function $\zeta^{MT}(z)$. Because it is not easy to establish a direct
relation between $D(z)$ and $\zeta^L(z)$, we introduce another
determinant, $L(z)$, called the homological determinant of $F$. These two
determinants are defined in a very similar way following the techniques
introduced in \cite{ASR}. To any $F$ we associate two pairs of linear endomorphisms 
$(\epsilon F_{\#0},\epsilon F_{\#1})$ and $(F_{\#0},F_{\#1})$. Although these
endomorphisms have in general infinite rank, we prove that their difference
has always finite rank. This allows us to define $D(z)$ and $L(z)$ as the
determinants of these pairs of linear endomorphisms.

This approach is different from the ones used by Milnor and Thurston (see 
\cite{MT}), Baladi and Ruelle (see \cite{BR}) or Baillif (see \cite{B}). For
better readability, we present basic algebraic notions and constructions in
the Appendix.

In the remainder of the paper we use the symbol $\Omega$ to denote a
finite and disjoint union of compact intervals on the real line 
\begin{equation*}
\Omega =[a_1,b_1]\cup [a_2,b_2]\cup \dots \cup [a_m,b_m]
\end{equation*}
with $a_1<b_1<a_2<b_2<\dots <a_m<b_m$.

\begin{definition}
A {\sf\em piecewise monotone} (shortly {\sf\em PM}) map on $\Omega$ 
is a map $F: \Omega \setminus C_F \to \Omega$ where $C_F$ is a
finite subset of $\Omega$ containing $\partial \Omega = \{a_1,b_1,\dots,a_m,b_m\}$ 
and such that $F$ is continuous and strictly
monotone on each connected component of $\Omega \setminus C_F$.
\end{definition}

Let $F: \Omega \setminus C_F \to \Omega$ be a PM map and $I=[x,y]$
(with $x<y$) be an interval. We say that $F$ is monotone on $I$ if
$]x,y[\subseteq \Omega \setminus C_F$. In this case we
define the sign function $\epsilon ([x,y])=\pm 1$ according to whether $F$
is increasing or decreasing on $]x,y[$. Moreover, for any 
$x \in \Omega \setminus C_F$, put $\epsilon (x)=\pm 1$ according to whether $F$ is
increasing or decreasing on a neighborhood of $x$ and put $\epsilon (x) = 0$
for every $x \in C_F$. By definition, a {\sf\em lap} of $F$ is a maximal
interval of monotonicity of $F$. That is to say, an interval $I = [c,d] \subseteq \Omega$ 
(with $c<d$) is a lap of $F$ if and only if $[c,d]\cap C_F=\{ c,d\}$. 
In what follows we use the symbol $\mathcal{L}_F$ to denote the set of all laps of $F$.

For a PM map $F:\Omega \setminus C_F \to \Omega$ and $n$ a
positive integer we define its $n$th iterate as a map 
$F^n: \Omega \setminus C_{F^n}\to \Omega$ inductively by 
$F^n(x)=F(F^{n-1}(x))$ for all $x\in \Omega \setminus C_{F^n}$ where 
\begin{equation*}
C_{F^n} = \{ x\in \Omega : F^k(x) \in C_F \text{ for some }k=0,\dots,n-1\} .
\end{equation*}
It can be easily seen that this map is PM as well.

Since it is easier to work with PM maps on the real line than with PMG maps
we want to replace the latter ones by the first ones. In fact, every PMG map
is induced by some PM map on an appropriate set $\Omega$ in the sense of the
following

\begin{definition}
\label{d4} Let $f:G\to G$ be a PMG map, $F:\Omega \setminus
C_F\to \Omega $ a PM map and $\pi :\Omega \to G$ a
continuous map such that $\mathrm{Ver}(G)\subseteq \pi (\partial \Omega )$
and $\pi $ maps $\Omega \setminus \partial \Omega $ homeomorphicaly into $
G\setminus \pi (\partial \Omega )$. Then we say that $f$ is {\sf\em{
induced}} by $F$ if the following diagram 
\begin{equation*}
\begin{CD}
\Omega \setminus C_F @>F>> \Omega\\ 
@V{\pi}VV @VV{\pi}V \\ 
G @>f>> G
\end{CD}
\end{equation*}
commutes.
\end{definition}

\subsection{The determinants $D(z)$ and $L(z)$}

Let $X$ be a topological space. Denote by $S_0(X;\mathbb{R})$ the 
$\mathbb{R}$-vector space whose basis consists of the formal symbols $x \in X$, and by
$S_1(X)$ its subspace generated by the vectors $y-x$ where $x$ and $y$ are
points lying in the same connected component of $X$. If $Y$ is a subset of $X$
and $F: X \setminus Y \to X$ is a map, we denote by 
$F_{\#0}: S_0(X;\mathbb{R}) \to S_0(X;\mathbb{R})$ the unique
linear endomorphism veryfying: $F_{\#0}(x)=F(x)$ if $x\in X \setminus Y$,
and $F_{\#0}(x)=0$ if $x\in Y$.

Let $F:\Omega \setminus C_F\to \Omega $ be a PM map. According to
the previous definitions, we have then a vector space $S_0(\Omega ;\mathbb{R})$, 
a subspace $S_1(\Omega ;\mathbb{R})$ of $S_0(\Omega ;\mathbb{R})$,
and a linear endomorphism $F_{\#0}:S_0(\Omega ;\mathbb{R})\to
S_0(\Omega ;\mathbb{R})$. Notice that both spaces $S_0(\Omega ;\mathbb{R}
)$ and $S_1(\Omega ;\mathbb{R})$ are infinite-dimensional but the quotient
space $S_0(\Omega ;\mathbb{R})/S_1(\Omega ;\mathbb{R})$ is
finite-dimensional with the dimension equal to the number of connected components
of $\Omega$.

Starting from $F_{\#0}$ we define another linear endomorphism 
$\epsilon F_{\#0}:S_0(\Omega ;\mathbb{R})\to S_0(\Omega ;\mathbb{R})$,
putting $\epsilon F_{\#0}(x)=\epsilon _F(x)F_{\#0}(x)$ for all $x\in \Omega$.
Next we define the linear endomorphisms 
$F_{\#1}:S_1(\Omega ;\mathbb{R})\to S_1(\Omega ;\mathbb{R})$ and 
$\epsilon F_{\#1}:S_1(\Omega ;\mathbb{R})\to S_1(\Omega ;\mathbb{R})$. Notice that, since $F$ is
a PM map, the subset of $S_1(\Omega ;\mathbb{R})$ 
\begin{equation*}
\mathcal I_F = \{ y-x:F\text{ is monotone on }[x,y]\}
\end{equation*}
spans $S_1(\Omega ;\mathbb{R})$. Furthermore if $F$ is monotone on $[x,y]$
then $F(y-)$ and $F(x+)$ lie in the same connected component of $\Omega $
and therefore $(F(y-)-F(x+))\in S_1(\Omega ;\mathbb{R})$ where $F(y-)$ and
$F(x+)$ denote the corresponding one-sided limits. So we can define $F_{\#1}$
and $\epsilon F_{\#1}$ as the unique linear endomorphisms of 
$S_1(\Omega ;\mathbb{R})$ such that
\begin{equation}\label{eqj0060}
F_{\#1}(y-x)=F(y-)-F(x+) \quad \text{and} \quad \epsilon F_{\#1}(y-x)=\epsilon
_F([x,y])F_{\#1}(y-x)  
\end{equation}
for all $y-x \in \mathcal I_F$.

As mentioned above, if $F:\Omega \setminus C_F\to \Omega $ is a PM map
then $F^n: \Omega \setminus C_{F^n} \to \Omega$ is also a PM map and
therefore the linear endomorphisms $F_{\#0}^n$, $F_{\#1}^n$, 
$\epsilon F_{\#0}^n$ and $\epsilon F_{\#1}^n$ are defined as well. The next lemma is a
simple consequence of the definitions and shows that the correspondences
$(\cdot)_{\#0}$, $(\cdot)_{\#1}$, $\epsilon(\cdot)_{\#0}$ and $\epsilon(\cdot)_{\#1}$
behave nicely under iteration.

\begin{lemma}\label{p5}
Let $F: \Omega \setminus C_F\to \Omega$ be a PM map. Then we
have $( F_{\#0})^n=F_{\#0}^n$, $(F_{\#1})^n=F_{\#1}^n$, 
$(\epsilon F_{\#0})^n=\epsilon F_{\#0}^n$
and $(\epsilon F_{\#1})^n=\epsilon F_{\#1}^n$ for all $n \geq 1$.
\end{lemma}

Thus for each PM map $F:\Omega \setminus C_F \to \Omega$ we have
two pairs of linear endomorphisms on $S_0(\Omega ;\mathbb{R})$, $(F_{\#0},F_{\#1})$ and 
$(\epsilon F_{\#0},\epsilon F_{\#1})$. Next we prove
that these pairs have both finite ranks (see Definition \ref{ad1}). For this
we need first to define extensions of $F_{\#1}$ and $\epsilon F_{\#1}$ to
the common superspace $S_0(\Omega ;\mathbb{R})$.

For each $c \in [a_i,b_i] \subset \Omega$, let 
$\alpha_c^-: \Omega \to \mathbb{R}$ and 
$\alpha_c^+: \Omega \to \mathbb{R}$ be step functions defined by
\begin{equation*}
\alpha_c^-(x)=
\begin{cases}
1 &\text{ for }x \in [c,b_i] \\ 0 &\text{ otherwise}
\end{cases}
\quad \text{and} \quad
\alpha_c^+(x)=
\begin{cases}
-1 &\text{ for }x \in ]c,b_i] \\ 0 &\text{ otherwise}.
\end{cases}
\end{equation*}
These step functions induce the linear forms
$\omega_c^-:S_0(\Omega ;\mathbb{R})\to \mathbb{R}$ and 
$\omega_c^+:S_0(\Omega ;\mathbb{R})\to \mathbb{R}$ defined by
$\omega_c^-(x)=\alpha_c^-(x)$ and $\omega_c^+(x)=\alpha_c^+(x)$
for all $x \in \Omega $. We also introduce the following
notation for special vectors from $S_0(\Omega ;\mathbb{R})$:
$v_c^-=F(c-)$, $\epsilon v_c^-=\epsilon _F(c-)v_c^-$,
$v_c^+=F(c+)$ and $\epsilon v_c^+=\epsilon _F(c+)v_c^+$ putting 
$v_c^-=0$ if $c=a_i$, and $v_c^+=0$ if $c=b_i$.

\begin{lemma}\label{p1}
Let $F$ be a PM map on $\Omega$. Then the linear endomorphism
$\varphi :S_0(\Omega ;\mathbb{R})\to S_0(\Omega ;\mathbb{R})$
defined by 
$\varphi =F_{\#0}+\sum_{c\in C_F}\omega _c^-\otimes v_c^-+\omega _c^+\otimes v_c^+$
is an extension of $F_{\#1}$to $S_0(\Omega ;\mathbb{R})$ that verifies
$\varphi (S_0(\Omega ;\mathbb{R}))\subseteq S_1(\Omega ;\mathbb{R})$. 
Consequently, the pair of linear
endomorphisms $(F_{\#0},F_{\#1})$ has finite rank and 
\begin{equation*}
\tr (F_{\#0},F_{\#1})=\sum_{c\in C_F}\omega _c^-\left(
v_c^-\right) +\omega _c^+\left( v_c^+\right).
\end{equation*}
\end{lemma}

\begin{proof}
Let $\varphi :S_0(\Omega ;\mathbb{R})\to S_0(\Omega ;\mathbb{R})$
be the linear endomorphism defined by $\varphi (x)=F_{\#1}(x-a_i)$, for
all $x\in \lbrack a_i,b_i]$. Evidently, $\varphi $ is an extension of $
F_{\#1}$to $S_0(\Omega ;\mathbb{R})$. Furthermore, since $\varphi
(a_i)=0 $ for all $i=1,\dots ,m$, it follows $\varphi \left( S_0(\Omega ;
\mathbb{R})\right) \subseteq S_1(\Omega ;\mathbb{R})$. On the other hand,
from the definitions of $F_{\#0}$ and $F_{\#1}$, we may write
\begin{equation*}
\begin{split}
F_{\#1}(x-a_i) &= F_{\#0}(x)-F_{\#0}(a_i)+\sum_{c\in ]a_i,x] \cap C_F}
v_c^- - \sum_{c \in [a_i,x[ \cap C_F} v_c^+ \\
&= F_{\#0}(x)+\sum_{c\in C_F}\omega _c^-(x) \ v_c^- + \omega _c^+ (x)  \ v_c^+
\end{split}
\end{equation*}
for all $x \in [a_i,b_i]$ and thus
\begin{equation*}
\varphi (x)=F_{\#0}(x)+\sum_{c\in C_F}\omega _c^-( x) \ v_c^-
+ \omega _c^+ (x) \ v_c^+
\end{equation*}
for all $x \in \Omega$.
\end{proof}

In the same way we can prove the following

\begin{lemma}\label{p2}
Let $F$ be a PM map on $\Omega $. Then the linear endomorphism
$\psi: S_0(\Omega; \mathbb{R}) \to S_0(\Omega; \mathbb{R})$
defined by $\psi = \epsilon F_{\#0} + \sum_{c \in C_F} \epsilon
\omega_c^- \otimes v_c^- + \epsilon \omega _c^+ \otimes v_c^+$ is an
extension of $\epsilon F_{\#1}$ to $S_0(\Omega; \mathbb{R})$ that verifies
$\psi (S_0(\Omega; \mathbb{R})) \subseteq S_1(\Omega; \mathbb{R})$.
Consequently, the pair of linear endomorphisms $(\epsilon F_{\#0},\epsilon F_{\#1})$
has finite rank and 
\begin{equation*}
\tr (\epsilon F_{\#0},\epsilon F_{\#1})=\sum_{c\in C_F}\omega_c^-(\epsilon v_c^-)
+ \omega_c^+(\epsilon v_c^+).
\end{equation*}
\end{lemma}

The last lemma shows that the determinants of the pairs $(\epsilon
F_{\#0},\epsilon F_{\#1})$ and $(F_{\#0},F_{\#1})$ (see the Appendix) are
defined. We define the {\sf\em{kneading determinant}} of $F$, $D(z)$,
and the {\sf\em{homological determinant}} of $F$, $L(z)$, by 
\begin{equation}\label{eqj0040}
\begin{split}
D(z)& =D_{(\epsilon F_{\#0},\epsilon F_{\#1})}(z) \\
& =\exp -\sum_{n\geq 1}\tr ((\epsilon F_{\#0})^n,(\epsilon F_{\#1})^n)\frac{z^n}n
\end{split}
\end{equation}
and
\begin{equation*}
\begin{split}
L(z)& =D_{(F_{\#0},F_{\#1})}(z) \\
& =\exp -\sum_{n\geq 1}\tr ((F_{\#0})^n,(F_{\#1})^n)\frac{z^n}n
\end{split}
\end{equation*}
Due to Lemmas \ref{p1}, \ref{p2} and Proposition \ref{ap4} there are vectors 
$u_1,\dots ,u_{p}\in S_0(\Omega ;\mathbb{R})$\textbf{\ }and linear forms 
$\nu _1$,\dots ,$\nu _{p}$, $\mu _1$,\dots , $\mu _{p}\in S_0(\Omega ;
\mathbb{R})^{\ast }$ such that 
\begin{equation*}
D(z)=\det (\mathbf{Id}-z\mathbf{M}(z))\text{ and }L(z)=\det (\mathbf{Id}-z
\mathbf{N}(z)),
\end{equation*}
where $\mathbf{M}(z)=[\mathbf{m}_{ij}(z)]$ and $\mathbf{N}(z)=[\mathbf{n}_{ij}(z)]$
are $p \times p$ matrices with entries from
$\mathbb{Z}[[z]]$ defined by 
\begin{equation}\label{eqj0065}
\mathbf{m}_{ij}(z)=\sum_{n\geq 0}\nu _i\circ (\epsilon F_{\#0})^n\left(
u_{j}\right) z^n \quad \text{and} \quad \mathbf n_{ij}(z)=\sum_{n\geq 0}\mu
_i\circ (F_{\#0})^n\left( u_{j}\right) z^n.
\end{equation}
Remark that as a consequence of the definitions the entries of
$\mathbf{M}(z)$ and $\mathbf{N}(z)$ are formal power series
that can be computed in terms of the orbits of the points of $C_F$ and
whose coefficients are from $\{-1,0,1\}$. Therefore the entries of
$\mathbf{M}(z)$ and $\mathbf{N}(z)$ and the corresponding
determinants $D(z)$ and $L(z)$ converge for all $|z|<1$.

At first glance, it is not clear which kind of relationship can hold between
the traces of $(F_{\#0},F_{\#1})$ and $(\epsilon F_{\#0},\epsilon F_{\#1})$
and the number of fixed points of $F$. For convenience we introduce the following notation.
Let the symbols $\mathcal{L}_F^+$ and $\mathcal{L}_F^-$
denote the set of all laps on which $F$ is increasing and
decreasing, respectively. We have then
$\mathcal{L}_F=\mathcal{L}_F^- \cup \mathcal{L}_F^+$. For each $I=[c,d] \in \mathcal{L}_F$,
define the number
\begin{equation*}
\sigma(I) = \omega_c^ + \left( v_c^+\right) + \omega_d^-\left(v_d^-\right).
\end{equation*}
Notice that from Lemmas \ref{p1} and \ref{p2} we have 
\begin{equation}\label{eqj0052}
\tr(F_{\#0},F_{\#1}) = \sum_{I \in \mathcal{L}_F} \sigma(I)
\end{equation}
and 
\begin{equation}\label{eqj0053}
\tr(\epsilon F_{\#0},\epsilon F_{\#1}) = \sum_{I \in \mathcal{L}_F} \epsilon_F(I)\sigma (I).
\end{equation}
On the other hand, it is easy to check the following

\begin{lemma}\label{p8}
Let $F: \Omega \setminus C_F\to \Omega$ be a PM map and 
$I=[c,d] \in \mathcal{L}_F$. Then $\sigma (I) \in \{ -1,0,1\}$
and
\begin{equation*}
\sigma (I) = 
\begin{cases} 1 &\text{if and only if $F(c+) \leq c$ and $d \leq F(d-)$}\\
-1 &\text{ if and only if $c<F(c+)$ and $F(d-)<d$}.
\end{cases}
\end{equation*}
\end{lemma}

We use this result to prove the following main relation between the
traces $\tr (\epsilon F_{\#0},\epsilon F_{\#1})$, and $\tr 
(F_{\#0},F_{\#1})$ and the number $\#\Fix ^-(F)$.

\begin{lemma}
\label{p6}Let $F:\Omega \setminus C_F\to \Omega $ be a PM map.
Then we have 
\begin{equation*}
\tr (\epsilon F_{\#0},\epsilon F_{\#1})-\tr 
(F_{\#0},F_{\#1})=2\#\Fix ^-(F)\text{.}
\end{equation*}
\end{lemma}

\begin{proof}
If $I=[c,d]\in \mathcal{L}_F^-$ then there is at most one fixed point of $F$ lying
in $\left] c,d\right[$ because $F$ is decreasing on $I$,
and from Lemma \ref{p8} we have: $\sigma (I)=-1$ if there exists a such
fixed point; $\sigma (I)=0$ otherwise. Therefore, from (\ref{eqj0052}) and (
\ref{eqj0053}), we obtain 
\begin{equation*}
\begin{split}
\tr (\epsilon F_{\#0},\epsilon F_{\#1})-\tr (F_{\#0},F_{\#1})
&=-2\sum_{I\in \mathcal{L}_F^-}\sigma (I) \\
&=2\#\Fix ^-(F)
\end{split}
\end{equation*}
as desired.
\end{proof}

Notice that by Lemmas \ref{p5} and \ref{p6} we have 
\begin{equation*}
\begin{split}
\tr ((\epsilon F_{\#0})^n,(\epsilon F_{\#1})^n)-\tr 
((F_{\#0})^n,(F_{\#1})^n) &=\tr (\epsilon F_{\#0}^n,\epsilon
F_{\#1}^n)-\tr (F_{\#0}^n,F_{\#1}^n) \\
&=2\#\Fix ^-(F^n)
\end{split}
\end{equation*}
and this proves the main theorem of this subsection,

\begin{theorem}
\label{t9}Let $F:\Omega \setminus C_F\to \Omega $ be a PM map.
Then 
\begin{equation*}
\exp \sum_{n \geq 1}\frac{2\#\Fix ^-(F^n)}n
z^n=L(z)D(z)^{-1}
\end{equation*}
holds in $\mathbb{Z}[[z]]$.
\end{theorem}

\subsection{Zeta functions and determinants}

Let $f:G\to G$ be a PMG map induced by $F:\Omega \setminus
C_F\to \Omega $. We have then two zeta functions, $\zeta ^-(z)$
and $\zeta ^{L}(z)$, and two determinants, $D(z)$ and $L(z)$.
By $P$ we denote the union of all periodic orbits of $f$ that
intersect $\pi(C_F)$ which is always a finite set.
Notice that the numbers $\#\Fix ^-(F^n)$
and $\#\Fix ^-(f^n)$ do not need to coincide because there may
exist periodic orbits of $f$ which intersect simultaneously $\pi (C_F)$
and $\#\Fix ^-(f^n)$ for some $n\geq 1$. Nevertheless we have
\begin{equation}\label{eqj0064}
\#\Fix ^-(f^n)-\#\Fix ^-(F^n)=\#P\cap \Fix ^-(f^n)
\end{equation}
for all $n \geq 1$, and consequently
\begin{equation}
\max \left\{ 1,\limsup_{n\to \infty }\#\Fix 
^-(F^n)^{1/n}\right\} =\max \left\{ 1,\limsup_{n\to \infty }\#
\Fix ^-(f^n)^{1/n}\right\}  \label{eqj0063}
\end{equation}
As an immediate consequence of (\ref{eqj0064}) and Theorem \ref{t9} we also
have:

\begin{corollary}\label{c2}
Let $f:G\to G$ be a PMG map induced by $F:\Omega \setminus
C_F\to \Omega $. Then 
\begin{equation*}
\zeta ^-(z)=L(z)D(z)^{-1}\exp \sum_{n\geq 1}\frac{2\#P\cap \mathrm{Fix
}^-(f^n)}nz^n
\end{equation*}
holds in $\mathbb{Z}[[z]]$.
\end{corollary}

The next result, together with Corollary \ref{c2}, allow us to establish a
main relationship between $\zeta^{MT}(z)$ and $D(z)$.

\begin{theorem}\label{l2}
Let $f: G\to G$ be a PMG map induced by $F: \Omega \setminus C_F \to \Omega $.
Then 
\begin{equation*}
\zeta^L(z) = L(z) \exp \sum_{n \geq 1} \frac{\#P \cap \Fix(f^n)}n z^n
\end{equation*}
holds in $\mathbb Z[[z]]$.
\end{theorem}

In order to prove Theorem \ref{l2}, let us begin by defining an auxiliary
pair of linear endomorphisms on $S_0(G;\mathbb{R})$. Let $f:G\to G$
be a PMG map induced by $F:\Omega \setminus C_F\to \Omega $ and
consider the linear endomorphisms $\beta _0:S_0(G;\mathbb{R})\to
S_0(G;\mathbb{R})$ and $\beta _1:S_1(G;\mathbb{R})\to S_1(G;
\mathbb{R})$ defined by \linebreak
$\beta_0(x)=
\begin{cases}
f(x) &\text{ if }x \in G \setminus \pi(C_F) \\
0 &\text{ if }x \in \pi(C_F)
\end{cases}$
and $\beta_1(y-x)=f\left( y\right) -f\left( x\right)$
for all $x$ and $y$ lying in the same connected component of $G$. The next
lemma shows that $(\beta _0,\beta _1)$ has finite rank and relates the
determinants $D_{(\beta _0,\beta _1)}(z)$ and $D_{f_{\ast 0}}(z)$.

\begin{lemma}
\label{p7}Under the conditions of the previous theorem, the pair $(\beta
_0,\beta _1)$ of linear endomorphisms on $S_0(G;\mathbb{R})$ has
finite rank and
\begin{equation*}
D_{(\beta _0,\beta _1)}(z)D_{f_{\ast 0}}(z)=\exp \sum_{n\geq 1}-
\frac{\#P\cap \Fix(f^n)}nz^n
\end{equation*}
holds in $\mathbb{Z}[[z]]$.
\end{lemma}

\begin{proof}
Let $\beta :S_0\left( G;\mathbb{R}\right) \to S_0\left( G;
\mathbb{R}\right) $ be the unique linear endomorphism that verifies $\beta
(x)=f\left( x\right) $ for all $x\in G$. Evidently, $\beta $ is an
extension of $\beta _1$ to $S_0\left( G;\mathbb{R}\right) $. Since $
\beta \left( x\right) =\beta _0\left( x\right) $ for all $x\in G\setminus
\pi \left( C_F\right) $ and $\pi \left( C_F\right) $ is a finite set,
we see that $\beta -\beta _0$ has finite rank. This shows that $(\beta
_0,\beta _1)$ has finite rank. On the other hand, because
\begin{equation*}
S_0\left( G;\mathbb{R}\right) /S_1\left( G;\mathbb{R}\right) =H_0(G;
\mathbb{R})\text{,}
\end{equation*}
we have the commutative diagram with exact rows
\begin{equation*}
\begin{CD}
0 @>>> S_1(G;\mathbb{R}) @>>> S_0(G;\mathbb{R}) @>{\pi_0}>> H_0(G;\mathbb{R}) @>>> 0 \\
@. @V{\beta_1}VV @V{\beta}VV @V{f_{\ast 0}}VV \\
0 @>>> S_1(G;\mathbb{R}) @>>> S_0(G;\mathbb{R}) @>{\pi_0}>> H_0(G;\mathbb{R}) @>>> 0
\end{CD}
\end{equation*}
Thus from Definition \ref{ad1} we have 
\begin{equation*}
\tr (\beta _0^n,\beta _1^n)+\tr (f_{\ast 0}^n)=
\tr (\beta ^n-\beta _0^n)
\end{equation*}
for all $n\geq 1$, and the proof follows because as an immediate consequence
of the definitions one has $\tr (\beta ^n-\beta _0^n)=\#P\cap \Fix(f^n)$ for all $n\geq 1$.
\end{proof}

Let us start now to prove Theorem \ref{l2}. Notice that if $f:G\to
G $ is a PMG map induced by $F:\Omega \setminus C_F\to \Omega$
then the continuous map $\pi :\Omega \to G$ induces linear
endomorphisms 
\begin{equation*}
\pi _0:S_0(\Omega ;\mathbb{R})\to S_0(G;\mathbb{R})\qquad 
\text{and}\qquad \pi _1:S_1(\Omega ;\mathbb{R})\to S_1(G;
\mathbb{R})
\end{equation*}
where $\pi _0$ is the unique linear map that verifies $\pi _0\left(
x\right) =\pi \left( x\right) $ for all $x\in \Omega $, and $\pi
_1$ is the restriction of $\pi _0$ to $S_1(\Omega ;\mathbb{R})$\textbf{
\ }(since $\pi $ is a continuous map, $\pi _0$ maps $S_1(\Omega ;\mathbb{
R})$\textbf{\ }in to $S_1(G;\mathbb{R})$). Since $\Ker \left( \pi
_1\right) =\Ker \left( \pi _0\right) \cap S_1(\Omega ;\mathbb{
R})\subset \Ker \left( \pi _0\right) $, we can consider the pair $
(\alpha _0,\alpha _1)$ of linear endomorphisms on $\Ker \left(
\pi _0\right) $, where $\alpha _i$ is the restriction of $F_{\#i}$ to $
\Ker \left( \pi _i\right) $. Notice that, since $\pi $ maps $
\Omega \setminus \partial \Omega $ homeomorphicaly into $G\setminus \pi
(\partial \Omega )$, we have $\Ker \left( \pi _0\right) \subset
S_0(\partial \Omega ;\mathbb{R})\subset \Ker \left(
F_{\#0}\right) $, and thus $\alpha _0=0$. On the other hand we also have $
H_1(G;\mathbb{R})=\Ker \left( \pi _1\right) $ and $\alpha
_1=f_{\ast 1}$. This shows that the pair $(\alpha _0,\alpha _1)$ has
finite rank and
\begin{equation}
D_{(\alpha _0,\alpha _1)}(z)=D_{(0,f_{\ast 1})}(z)=D_{f_{\ast 1}}(z)
\label{eqj0061}
\end{equation}

Thus we obtain three pairs of linear endomorphisms with finite rank, $(\alpha
_0,\alpha _1)$, $(F_{\#0},F_{\#1})$ and $(\beta _0,\beta _1)$ on $
\Ker \left( \pi _0\right) $, $S_0\left( \Omega ;\mathbb{R}
\right) $ and $S_0(G;\mathbb{R})$, respectively, and the two following
commutative diagrams with exact rows
\begin{equation}\label{eqj0043}
\begin{CD}
0 @>>> \Ker \pi_0 @>>> S_0(\Omega;\mathbb{R}) @>{\pi_0}>> S_0(G;\mathbb{R}) @>>> 0 \\
@. @V0VV @V{F_{\#0}}VV @V{f_{\#0}}VV \\
0 @>>> \Ker \pi_0 @>>> S_0(\Omega;\mathbb{R}) @>{\pi_0}>> S_0(G;\mathbb{R}) @>>> 0
\end{CD}
\end{equation}
and
\begin{equation}\label{eqj0044}
\begin{CD}
0 @>>> \Ker \pi_1 @>>> S_1(\Omega;\mathbb{R}) @>{\pi_1}>> S_1(G;\mathbb{R}) @>>> 0 \\
@. @V{\alpha_1}VV @V{F_{\#1}}VV @V{\beta_1}VV \\
0 @>>> \Ker \pi_1 @>>> S_1(\Omega;\mathbb{R}) @>{\pi_1}>> S_1(G;\mathbb{R}) @>>> 0
\end{CD}
\end{equation}
The restriction of $\pi _0$ to $S_1( \Omega ;\mathbb{R})$ is $\pi _1$, 
so from Proposition \ref{ap3} and Lemma \ref{p7} we obtain
\begin{equation*}
\begin{split}
L(z) &=D_{(F_{\#0},F_{\#1})}(z) \\
&=D_{(\alpha _0,\alpha _1)}(z)D_{(\beta _0,\beta _1)}(z) \\
&=D_{(\alpha _0,\alpha _1)}(z)D_{f_{\ast 0}}(z)^{-1}\exp \sum_{n\geq 1}-
\frac{\#P\cap \Fix(f^n)}n z^n
\end{split}
\end{equation*}
and from (\ref{eqj0060})
\begin{equation*}
\begin{split}
L(z) &=D_{f_{\ast 1}}(z)D_{f_{\ast 0}}(z)^{-1}\exp \sum_{n\geq 1}-\frac{
\#P\cap \Fix(f^n)}nz^n \\
&=\zeta ^{L}(z)\exp \sum_{n\geq 1}-\frac{\#P\cap \Fix(f^n)}n z^n
\end{split}
\end{equation*}
as desired.

From Corollary \ref{c2} and Theorem \ref{l2} we obtain

\begin{theorem}
\label{t10}Let $f:G\to G$ be a PMG map induced by $F:\Omega
\setminus C_F\to \Omega $. Then there exists a formal power series 
$H(z)$ such that $H(z)$ converges and is nonzero for all $\left\vert z\right\vert <1$ and 
\begin{equation*}
\zeta ^{MT}(z)=H(z)D(z)^{-1}
\end{equation*}
holds in $\mathbb{Z[[}z\mathbb{]]}$.
\end{theorem}

\begin{proof}
We have
\begin{equation*}
\zeta ^{MT}(z)=\zeta ^-(z)\zeta ^{L}(z)^{-1}=D(z)^{-1}\exp a(z)
\end{equation*}
with
\begin{equation*}
a(z)=\sum_{n\geq 1}\frac{2\#P\cap \Fix ^-(f^n)-\#P\cap \mathrm{Fix
}(f^n)}nz^n\text{.}
\end{equation*}
Thus, because $P$ is a finite set, it follows immediately that $a(z)$
converges for all $\left\vert z\right\vert <1$ and consequently $H(z)=\exp
a(z)\neq 0$ for all $\left\vert z\right\vert <1$.
\end{proof}

As mentioned before, the kneading determinant $D(t)$ converges for all $z\in 
\mathbb{D}=\left\{ z\in \mathbb C:\left\vert z\right\vert <1\right\} $, so,
as immediate consequence of Theorem \ref{t10}, we obtain

\begin{corollary}
\label{c3}Let $f:G\to G$ be a PMG map induced by $F:\Omega \setminus
C_F\to \Omega $, $\rho $ the radius of convergence of $\zeta
^{MT}(z)$, and $z_0$ a zero of $D(z)$ lying in $\mathbb{D}$. Then we have $
\rho \leq \left\vert z_0\right\vert $.
\end{corollary}

Let $f:G\to G$ be a PMG map induced by $F:\Omega \setminus
C_F\to \Omega $. Next we discuss the relationship between
topological entropy of $f$ and the zeros of $D\left( z\right) $. From now on
we use the symbol $\ell \left( F^n\right) $ to denote the number of
laps of the iterate $F^n:\Omega \setminus C_{F^n}\to \Omega $,
in other words $\ell \left( F^n\right) =\#\mathcal{L}(F^n)$. If $
I=[c,d]\in \mathcal{L}(F^n)$ we define the variation of $F^n$ on $I$ by $
\mathrm{Var}_i(F^n)=\left\vert F^n(d-)-F^n(d+)\right\vert $. The
variation of $F^n$ is defined by
\begin{equation*}
\mathrm{Var}(F^n)=\sum_{I\in \mathcal{L}(F^n)}\mathrm{Var}
_i(F^n)\text{.}
\end{equation*}
Recall that, for interval and circle maps Misiurewicz and Szlenk proved in
[MSz] that 
\begin{equation}
\htop (f)=\log \lim_{n\to \infty }\ell (F^n)^{1/n}=\log
\max \left\{ 1,\lim_{n\to \infty }\mathrm{Var}(F^n)^{1/n}\right\} 
\text{,}  \label{eqj0057}
\end{equation}
and the same arguments can be adapted to show that (\ref{eqj0057}) holds for
any PMG map.

Let us begin with the following

\begin{theorem}
\label{t8} Let $f:G\to G$ be a PMG map induced by $F:\Omega
\setminus C_F\to \Omega $. If $\htop (f)>0$. Then $
D(z)=0 $ for some $\left\vert z\right\vert =\lim_{n\to \infty }
\mathrm{Var}(F^n)^{-1/n}$.
\end{theorem}

\begin{proof}
Let $v=(b_1-a_1)+...+(b_{m}-a_{m})\in S_1(\Omega ;\mathbb{R})$ and $
\xi :S_1(\Omega ;\mathbb{R})\to \mathbb{R}$ the linear form
defined by $\xi (y-x)=y-x$, for all $x$ and $y$ lying in the same connected
component of $\Omega $. We have then a pair $\left( \epsilon
F_{\#0},\epsilon F_{\#1}+\xi \otimes v\right) $ of endomorphisms on $
S_1(\Omega ;\mathbb{R})$ with finite rank. As mentioned before, the
kneading determinant $D(z)=D_{\left( \epsilon F_{\#0},\epsilon
F_{\#1}\right) }(z)$ converges for all $\left\vert z\right\vert <1$. Using
the same argument, it is easy to show that $D_{\left( \epsilon
F_{\#0},\epsilon F_{\#1}+\xi \otimes v\right) }(z)$ also converges for all $
\left\vert z\right\vert <1$.

Notice that if $I=[c,d]$ is a lap of $F^n$, we have $\mathrm{Var}
_i(F^n)=\xi \circ \epsilon F_{\#1}^n(d-c)$. Thus, from the linearity
of $\epsilon F_{\#1}^n$, and by Lemma \ref{p5}, we arrive at 
\begin{equation*}
\mathrm{Var}(F^n)=\xi \circ \epsilon F_{\#1}^n(v)=\xi \circ \left(
\epsilon F_{\#1}\right) ^n(v)\text{,}
\end{equation*}
since $\htop (f)>0$, we have then 
\begin{equation*}
\underset{n\to \infty }{\lim \sup }\left\vert \xi \circ \left(
\epsilon F_{\#1}\right) ^n(v)\right\vert ^{1/n}=\lim_{n\to \infty }
\mathrm{Var}(F^n)^{1/n}>1\text{,}
\end{equation*}
and from Proposition \ref{ap5}, $D(z)=0$, for some $\left\vert z\right\vert
=\lim_{n\to \infty }\mathrm{Var}(F^n)^{-1/n}$.
\end{proof}

We have now everything we needed to prove Theorem \ref{t2}. Let $
f:G\to G$ be a PMG map induced by $F:\Omega \setminus
C_F\to \Omega $, and denote by $\rho $ the fradius of convergence
of $\zeta _F^{MT}(z)$. From the definition of $\zeta _F^{MT}(z)$, we see
at once that 
\begin{equation*}
\rho ^{-1}\leq \max \left\{ 1,\limsup_{n\to \infty }\#\mathrm{
Fix}^-(f^n)^{1/n},r(f_{\ast 1})\right\} \text{,}
\end{equation*}
and thus
\begin{equation*}
\log \max \left\{ 1,\rho ^{-1}\right\} \leq \max \left\{ h_{\mathrm{per}
}^-(f),h_{\mathrm{\hom }}(f)\right\} \text{.}
\end{equation*}
Suppose that $\htop (f)>0$. In this case, from (\ref{eqj0057}) and
Theorem \ref{t8}, we have $D(z)=0$, for some $\left\vert z\right\vert
=\lim_{n\to \infty }\mathrm{Var}(F^n)^{-1/n}<1$. Thus from
Corollary \ref{c3}, we have $\rho \leq \lim_{n\to \infty }\mathrm{Var
}(F^n)^{-1/n}$ and thus 
\begin{equation*}
\htop (f)\leq \log \max \left\{ 1,\rho ^{-1}\right\} \text{.}
\end{equation*}
Finally it remains to prove $\max \left\{ h_{\mathrm{per}}^-(f),h_{\mathrm{
\hom }}(f)\right\} \leq \htop (f)$. Since in each lap of $
F^n$ there is at most one fixed point of negative type, we have 
\begin{equation*}
\lim_{n\to \infty }\ell \left( F^n\right) ^{1/n}\geq \max \left\{
1,\underset{n\to \infty }{\lim \sup }\#\Fix 
^-(F^n)^{1/n}\right\} \text{,}
\end{equation*}
and from (\ref{eqj0063})
\begin{equation*}
\lim_{n\to \infty }\ell \left( F^n\right) ^{1/n}\geq \max \left\{
1,\underset{n\to \infty }{\lim \sup }\#\Fix 
^-(f^n)^{1/n}\right\}.
\end{equation*}
Thus, by (\ref{eqj0057}), $\htop (f)\geq h_{\mathrm{per}}^-(f)$.
From (\ref{eqj0052}) and Lemma \ref{p5} we also have 
\begin{equation*}
\begin{split}
\ell \left( F^n\right) &\geq \left\vert \sum_{I\in \mathcal{L}_{F^n}}\sigma (I)\right\vert \\
&=\left\vert \tr (F_{\#0}^n,F_{\#1}^n)\right\vert \\
&=\left\vert \tr (\left( F_{\#0}\right) ^n,\left( F_{\#1}\right)
^n)\right\vert
\end{split}
\end{equation*}
for all $n \geq 1$. But from Theorem \ref{l2}
\begin{equation*}
\left\vert \tr (\left( F_{\#0}\right) ^n,\left( F_{\#1}\right)
^n)\right\vert =\left\vert \#P\cap \Fix(f^n)+\tr (f_{\ast
1})^n-\tr (f_{\ast 0})^n\right\vert
\end{equation*}
for all $n\geq 1$. Thus, since $P$ is a finite set, it follows
\begin{equation*}
\lim_{n\to \infty }\ell \left( F^n\right) ^{1/n}\geq \max \left\{
1,\underset{n\to \infty }{\lim \sup }\left\vert \tr (f_{\ast
1})^n\right\vert ^{1/n}\right\} =r(f)\text{,}
\end{equation*}
and, once again from (\ref{eqj0057}), $\htop(f) \geq \hhom(f)$.

\section*{Appendix (Pairs of linear endomorphisms)}

Let $V$ be a vector space over $\mathbb{R}$ and let $\varphi: V \to V$
be a linear map with finite rank. As usually we define the trace of $\varphi$ by 
\begin{equation*}
\tr\varphi = \tr\varphi_{|\Image\varphi}. 
\end{equation*}
If $\varphi$ has finite rank then there are vectors $v_1,\dots,v_k\in V$
and linear forms $\omega_1,\dots,\omega_k \in V^\ast$ such that 
\begin{equation*}
\varphi = \sum_{i=1}^k\omega_i \otimes v_i. 
\end{equation*}
Considering the matrix
\begin{equation}\label{af1}
\mathbf{M}=
\begin{pmatrix}
\omega_1(v_1) & \dots & \omega_1(v_k) \\ 
\vdots &  & \vdots \\ 
\omega_k(v_1) & \dots & \omega_k(v_k)
\end{pmatrix}
\end{equation}
we have
\begin{equation*}
\tr\varphi =\tr\mathbf{M}. 
\end{equation*}
More generally, if $\varphi$ has finite rank then $\varphi^n$, $n\geq 1$,
has also finite rank and 
\begin{equation*}
\tr\varphi ^n=\tr\mathbf{M}^n. 
\end{equation*}
The following proposition is well known and gives an explicit method for
computing the numbers $\tr\varphi^n$ for $n \geq 1$. Defining the
determinant of $\varphi$ to be the following formal power series 
\begin{equation*}
D_\varphi(z) = \exp\sum_{n\geq 1} - \tr\varphi^n \frac{z^n}n
\end{equation*}
we have

\begin{proposition}
\label{ap1} Let $\varphi $ be an endomorphism with finite rank. Then we have 
\begin{equation*}
D_\varphi(z) = \det(\mathbf{Id}-z\mathbf{M}) \qquad\text{in }\mathbb{R}[[z]].
\end{equation*}
\end{proposition}

Now we consider a more general situation. By a pair of endomorphisms
$(\varphi_0,\varphi_1)$ on $V$ we mean two two linear maps
$\varphi_0: V_0 \to V_0$ and $\varphi_1: V_1 \to V_1$ defined on two
finite-codimensional subspaces $V_0$ and $V_1$ of the same $\mathbb{R}$-vector space $V$.

\begin{definition}
\label{ad1} We say that the pair of endomorphisms $(\varphi _0,\varphi
_1)$ on $V$ has a {\sf\em finite rank} if there exist linear maps $\overline{
\varphi }_0$, $\widetilde{\varphi }_0$, $\overline{\varphi }_1$ and $
\widetilde{\varphi }_1$ such that the following diagram with exact rows 
\begin{equation*}
\begin{CD}
0 @>>> V_j @>\subseteq>> V @>{\mathrm{pr}}>> V/V_j @>>> 0 \\
@. @VV{\varphi_j}V @VV{\overline{\varphi}_j}V @VV{\widetilde{\varphi}_j}V \\
0 @>>> V_j @>\subseteq>> V @>{\mathrm{pr}}>> V/V_j @>>> 0
\end{CD}
\end{equation*}
commutes for $j=0,1$ and the linear map $\overline{\varphi}_1-\overline{\varphi}_0$
has finite rank. The {\sf\em trace} of a pair $(\varphi_0,\varphi_1)$
with finite rank is defined by $\tr(\varphi_0,\varphi_1)=
\tr(\overline{\varphi}_1-\overline{\varphi}_0)-\tr\widetilde{\varphi}_1+\tr\widetilde{\varphi}_0$.
\end{definition}

It is easy to see that the definition does not depend on $\overline{\varphi}_0$,
$\widetilde{\varphi}_0$, $\overline{\varphi}_1$ and $\widetilde{\varphi}_1$.
As an immediate consequence of the definition we get

\begin{proposition}
\label{ap2} Let $(\varphi _{0},\varphi _{1})$ be a pair of endomorphisms on $
V$, and $\overline{\varphi }_{i}:V\rightarrow V$ an extension of $\varphi
_{i}$ such that $\overline{\varphi }_{i}\left( V\right) \subseteq V_{i}$,
for $i=0,1$. Then $(\varphi _{0},\varphi _{1})$ has finite rank if and only
if $\overline{\varphi }_{1}-\overline{\varphi }_{0}$ has finite rank.
Furthermore, if $(\varphi _{0},\varphi _{1})$ has finite rank then $\mathrm{
tr}(\varphi _{0},\varphi _{1})=\mathrm{tr}(\overline{\varphi }_{1}-\overline{
\varphi }_{0})$.
\end{proposition}

Let $(\varphi _0,\varphi _1)$ be a pair of endomorphisms on $V$ having
finite rank, and consider endomorphisms $\overline{\varphi }_0$ and $
\overline{\varphi }_1$ as in Proposition \ref{ap2}. Since $\overline{
\varphi }_1-\overline{\varphi }_0$ has finite rank, we can consider
vectors $v_1,\dots ,v_{k}\in V$ and linear forms $\omega _1,\dots
,\omega _{k}\in V^{\ast }$ such that 
\begin{equation}
\overline{\varphi }_1-\overline{\varphi }_0=\sum_{i=1}^{k}\omega
_i\otimes v_i.  \label{af2}
\end{equation}
More generally, we have 
\begin{equation*}
\overline{\varphi }_1^n-\overline{\varphi }_0^n=\sum_{i=1}^{k}
\sum_{j=1}^n\left( \omega _i\circ \overline{\varphi }_1^{n-j}\right)
\otimes \overline{\varphi }_0^{j-1}(v_i), 
\end{equation*}
for each $n\geq 1$. This shows that $\overline{\varphi }_1^n-\overline{
\varphi }_0^n$ has finite rank for each $n\geq 1$. Thus, once more from
Proposition \ref{ap2}, we conclude that the pair $(\varphi _0^n,\varphi
_1^n)$ has finite rank and 
\begin{equation*}
\tr(\varphi _0^n,\varphi _1^n)=\tr(\overline{\varphi }_1^n-
\overline{\varphi }_0^n)\qquad \text{for each }n\geq 1. 
\end{equation*}

\begin{definition}
\label{ad2} Let $(\varphi_0, \varphi_1)$ be a pair of endomorphisms having
finite rank. We define the {\sf\em determinant} of $(\varphi_0, \varphi_1)$ to
be the following element of $\mathbb{R}[[z]]$ 
\begin{equation*}
D_{(\varphi_0, \varphi_1)}(z) = - \sum_{n \geq 1} \tr (\varphi_0^n,
\varphi_1^n) \frac{z^n}n. 
\end{equation*}
\end{definition}

Observe that if $\varphi$ has finite rank then 
\begin{equation*}
D_{(0, \varphi)}(z) = D_\varphi(z). 
\end{equation*}
If $\varphi_0$ and $\varphi_1$ both have finite ranks then 
\begin{equation*}
D_{(\varphi_0, \varphi_1)}(z) = D_{\varphi_1}(z) D_{\varphi_0}(z)^{-1}. 
\end{equation*}
So, in these cases, we can use Proposition \ref{ap1} for computing $
D_{(\varphi_0, \varphi_1)}(z)$. Obviously, in the general case, Proposition 
\ref{ap1} does not allow us to compute $D_{(\varphi_0, \varphi_1)}(z)$ --- $
D_{\varphi_0}(z)$ and $D_{\varphi_1}(z)$ are not defined in general.

In order to compute $D_{(\varphi_0, \varphi_1)}(z)$ in the general case, we
generalize the Proposition \ref{ap1}. Let $\overline{\varphi}_0$ and $
\overline{\varphi}_1$ be endomorphisms as in Proposition \ref{ap2}.
Considering vectors $v_1,\dots,v_k \in V$ and linear forms $
\omega_1,\dots,\omega_k \in V^\ast$ as in (\ref{af2}), we define the matrix

\begin{equation}
\mathbf{M}(z) = 
\begin{pmatrix}
\sum_{n \geq 0} \omega_1 (\overline{\varphi}_0^n(v_1)) z^n & \dots & \sum_{n
\geq 0} \omega_1 (\overline{\varphi}_0^n(v_k) z^n \\ 
\vdots &  & \vdots \\ 
\sum_{n \geq 0} \omega_k (\overline{\varphi}_0^n(v_1)) z^n & \dots & \sum_{n
\geq 0} \omega_k (\overline{\varphi}_0^n(v_k) z^n \label{af3}
\end{pmatrix}
\end{equation}
with coefficients in $\mathbb{R}[[z]]$. Observe that if we identify an
endomorphism with finite rank $\varphi: V \to V$ with the corresponding pair
of finite rank $(0, \varphi)$ then the matrix $\mathbf{M}(z)$ from (\ref{af3}
) coincides with the matrix $\mathbf{M}$ defined in (\ref{af1}). Thus the
next proposition which gives an explicit method for computing $D_{(\varphi,
\psi)}(z)$, can be regarded as a generalization of Proposition \ref{ap1}.

\begin{proposition}
\label{ap4} Let $(\varphi _0,\varphi _1)$ be a pair of endomorphisms
having finite rank. Then
\begin{equation*}
D_{(\varphi_0,\varphi_1)}(z) = \det(\mathbf{Id}-z\mathbf{M}(z)).
\end{equation*}
\end{proposition}

Let $0\to U\overset{i}{\to }V\overset{p}{\to }
W\to 0$ be an exact sequence of $\mathbb{R}$-vector spaces. Recall
that if the diagram 
\begin{equation*}
\begin{CD}
0 @>>> U @>i>> V @>p>> W @>>> 0 \\
@. @VV{\chi}V @VV{\varphi}V @VV{\psi}V \\
0 @>>> U @>i>> V @>p>> W @>>> 0
\end{CD}
\end{equation*}
commutes and the endomorphisms $\chi $, $\varphi $ and $\psi $ have finite
rank then we have 
\begin{equation*}
\tr\varphi ^n=\tr\chi ^n+\tr\psi ^n 
\end{equation*}
for all $n\geq 1$, and therefore 
\begin{equation*}
D_{\varphi }(z)=D_{\chi }(z)D_{\psi }(z)\text{.} 
\end{equation*}
The next proposition can be regarded as a generalization of the last
formula. Let $(\chi _0,\chi _1)$, $(\varphi _0,\varphi _1)$ and $
(\psi _0,\psi _1)$ be pairs of endomorphisms in $U$, $V$ and $W$
respectively, such that the following diagram 
\begin{equation*}
\begin{CD}
0 @>>> U_j @>i>> V_j @>p>> W_j @>>> 0 \\
@. @VV{\chi_j}V @VV{\varphi_j}V @VV{\psi_j}V \\
0 @>>> U_j @>i>> V_j @>p>> W_j @>>> 0
\end{CD}
\end{equation*}
commutes for $j=0,1$. Then we have the following

\begin{proposition}\label{ap3}
Let $(\chi_0,\chi_1)$, $(\varphi_0,\varphi_1)$ and $(\psi_0,\psi_1)$
be pairs of endomorphisms with finite rank such that
the last diagram above commutes. Then
\begin{equation*}
\begin{split}
\tr(\varphi_0^n,\varphi_1^n) &= \tr(\chi_0^n,\chi_1^n)+\tr(\psi_0^n,\psi_1^n)
\qquad \text{for all }n\geq 1,\text{ and}\\
D_{(\varphi_0,\varphi_1)}(z) &= D_{(\chi_0,\chi_1)}(z)D_{(\psi_0,\psi_1)}(z).
\end{split}
\end{equation*}
\end{proposition}

Let us consider, for the last time, a linear endomorphism $\varphi
:V\to V$ with finite rank. Recall that, if $v\in V$ and $\xi \in
V^{\ast }$, then there exists an eigenvalue, $\lambda $, of $\varphi $ such
that
\begin{equation*}
\limsup_{n\to \infty }\left\vert \xi \circ \varphi
^{n}(v)\right\vert ^{1/n}=\left\vert \lambda \right\vert
\end{equation*}
and consequently
\begin{equation}
D_{\varphi }(z)=0\text{ for some }\left\vert z\right\vert =\frac{1}{
\limsup_{n\to \infty }\left\vert \xi \circ \varphi
^{n}(v)\right\vert ^{1/n}}\text{. }  \label{eqj0066}
\end{equation}
We will finish this appendix with a generalization of (\ref{eqj0066}). Let $
(\varphi _{0},\varphi _{1})$ be a pair of endomorphisms on $V$ with finite
rank, $v\in V_{1}$, $\xi \in V_{1}^{\ast }$. Then the pair $(\varphi
_{0},\varphi _{1}+\xi \otimes v)$ of endomorphisms on $V$, also has finite
rank. Notice that, since $D_{(\varphi _{0},\varphi _{1})}(z)$ and $
D_{(\varphi _{0},\varphi _{1}+\xi \otimes v)}(z)$ are not (in general)
polynomials, we have to assumme that there exists $r>0$ such that $
D_{(\varphi _{0},\varphi _{1})}(z)$ and $D_{(\varphi _{0},\varphi _{1}+\xi
\otimes v)}(z)$ converge for all $\left\vert z\right\vert <r$. If we
consider the pair $(\varphi _{1},\varphi _{1}+\xi \otimes v)$ of
endomorphisms on $V_{1}$, this pair has evidently finite rank, and from
Proposition \ref{ap4} we see that
\begin{equation*}
D_{(\varphi _{1},\varphi _{1}+\xi \otimes v)}(z)=1-\sum_{n\geq 0}\xi \circ
\varphi _{1}^{n}(v)z^{n+1}\text{ }
\end{equation*}
holds in $\mathbb{R}[[z]]$. On the other hand, regarding $(\varphi
_{1},\varphi _{1}+\xi \otimes v)$ as a a pair of endomorphisms on $V$, we
have the decomposition
\begin{equation*}
D_{(\varphi _{1},\varphi _{1}+\xi \otimes v)}(z)=D_{(\varphi _{1},\varphi
_{0})}(z)D_{(\varphi _{0},\varphi _{1}+\xi \otimes v)}(z)=\frac{D_{(\varphi
_{0},\varphi _{1}+\xi \otimes v)}(z)}{D_{(\varphi _{0},\varphi _{1})}(z)}
\text{,}
\end{equation*}
and therefore
\begin{equation*}
1-\sum_{n\geq 0}\xi \circ \varphi _{1}^{n}(v)z^{n+1}=\frac{D_{(\varphi
_{0},\varphi _{1}+\xi \otimes v)}(z)}{D_{(\varphi _{0},\varphi _{1})}(z)}
\end{equation*}
holds in in $\mathbb{R}[[z]]$. Thus, since the radius of convergence of 
\begin{equation*}
1-\sum_{n\geq 0}\xi \circ \varphi _{1}^{n}(v)z^{n+1}
\end{equation*}
is
\begin{equation*}
\rho =\frac{1}{\limsup_{n\to \infty }\left\vert \xi \circ \varphi
_{1}^{n}(v)\right\vert ^{1/n}}\text{,}
\end{equation*}
and the function 
\begin{equation*}
\gamma \left( z\right) =\frac{D_{(\varphi _{0},\varphi _{1}+\xi \otimes
v)}(z)}{D_{(\varphi _{0},\varphi _{1})}(z)}
\end{equation*}
is meromorphic on $\left\vert z\right\vert <r$, we can conclude: if $\rho <r$
then $\gamma \left( z\right) $ has a pole lying in $\left\vert z\right\vert
=\rho $. So, because the poles of $\gamma \left( z\right) $ are zeros of $
D_{(\varphi _{0},\varphi _{1})}(z)$, we may write:

\begin{proposition}
\label{ap5}Let $(\varphi _{0},\varphi _{1})$ be a pair of endomorphisms on $
V $ with finite rank, $v\in V_{1}$, $\xi \in V_{1}^{\ast }$ and $r>0$ such
that $D_{(\varphi _{0},\varphi _{1})}(z)$ and $D_{(\varphi _{0},\varphi
_{1}+\xi \otimes v)}(z)$ converge for all $\left\vert z\right\vert <r$, \\and $
\limsup_{n\to \infty }\left\vert \xi \circ \varphi
_{1}^{n}(v)\right\vert ^{1/n}>r^{-1}$. Then we have 
\begin{equation*}
D_{(\varphi _{0},\varphi _{1})}(z)=0\text{, for some }\left\vert
z\right\vert =\frac{1}{\limsup_{n\to \infty }\left\vert \xi \circ
\varphi _{1}^{n}(v)\right\vert ^{1/n}}
\end{equation*}
\bigskip
\end{proposition}

\end{document}